\documentclass[10pt]{article}
\usepackage{latexsym}
\usepackage{amsmath}
\renewcommand{\epsilon}{\varepsilon}

\renewcommand{\baselinestretch}{2}

\usepackage{amsmath,bm}
\usepackage{amssymb}
\usepackage[ansinew]{inputenc}
\usepackage{graphicx}
\usepackage{epsfig}
\usepackage{dsfont}
\usepackage{bbm}

\usepackage[authoryear]{natbib}
\newtheorem{satz}{Theorem}[section]

\newtheorem{rem}[satz]{Remark}

\def\3{\ss}

\def \R{I \!\! R}   
\def \N{I \!\! N}   

\newcommand{\E}{\mathbbm{E}}

\newcommand{\bea}{\begin{eqnarray*}}
\newcommand{\eea}{\end{eqnarray*}}
\newcommand{\be}{\begin{eqnarray}}
\newcommand{\ee}{\end{eqnarray}}

\newcommand{\ba}{\begin{array}}
\newcommand{\ea}{\end{array}}
\newcommand{\cum}{\text{cum}}

\def\3{\ss}

\begin{document}

\begin{center}
{\LARGE \textbf{Comparing spectral densities of stationary time series with unequal sample sizes}}
\setcounter{footnote}{0}

\bigskip
\Large
Philip Preu{\ss}\footnote{Address for correspondence: Philip Preu{\ss}, 
Ruhr-Universit\"at Bochum,
Fakult\"at f\"ur Mathematik, 
NA 3/30,
Universit\"atsstr. 150,
D-44780 Bochum, Germany,
email: philip.preuss@rub.de,
Fon: +49/234/32--23289,
Fax: +49/551/32--14559} and Thimo Hildebrandt\\*[0.2cm]
\large
Fakult\"at f\"ur Mathematik\\
Ruhr-Universit\"at Bochum, Germany\\*[0.2cm]
\today

\end{center}
\normalsize

\begin{abstract}
This paper deals with the comparison of several stationary processes with unequal sample sizes. We provide a detailed theoretical framework on the testing problem for equality of spectral densities in the bivariate case, after which the generalization of our approach to the $m$ dimensional case and to other statistical applications (like testing for zero correlation or clustering of time series data with different length) is straightforward. We prove asymptotic normality of an appropriately standardized version of the test statistic both under the null and the alternative and investigate the finite sample properties of our method in a simulation study. Furthermore we apply our approach to cluster financial time series data with different sample length.
\end{abstract}



\section{Introduction}
\def\theequation{1.\arabic{equation}}
\setcounter{equation}{0} The comparison and clustering of different time series is an important topic in statistical data analysis and has various applications in fields like economics, marketing, medicine and physics, among many others. Examples are the grouping of stocks in several categories for portfolio selection in finance or the identification of similar birth and death rates in population studies. One approach to identify similarities or dissimilarities between two stationary processes is to compare the spectral densities of both time series, which directly yields to the testing problem for equality of spectral densities in multivariate time series data. This problem has found considerable interest in the literature [see for example \cite{jenkins1961} or \cite{desouza1982} for some early results], but in the nonparametric situation nearly all proposed procedures are only reasoned by simulation studys or heuristic proofs, see \cite{coatesdiggle1986}, \cite{poetscherreschenhofer1988}, \cite{digglefischer1991} and \cite{maharaj2002} among many others. Most recently \cite{eichler2008}, \cite{detpap2009}, \cite{detkinsvet2010}, \cite{jentschpauly2011} and \cite{dethil2011} provided mathematical details for the above testing problem using different $L_2$-type statistics, but nevertheless in all mentioned articles it is always required that the different time series have the same length, which is typically not the case in practice. \cite{caiado2009} considered different metrics for the comparison of time series with unequal sample sizes in a simulation study and \cite{jentsch2011} provided a theoretical result, which however does not yield a consistent test as it was also pointed out by the authors. \\
This paper generalizes the approach of \cite{detkinsvet2010} to the case of unequal sample sizes and yields a consistent test for the equalness of spectral densities for time series with different length. Although the limiting distribution will be the same as in \cite{detkinsvet2010} note that our proof is completely different. This is due to the fact that one essential part in the proofs of \cite{detkinsvet2010} is that the different processes have the same Fourier coefficents which is not given if the observed time series have different sample sizes. For the sake of brevity we will focus on the case of two (not necessarily independent) stationary processes, but the results can be easily extended to the case of an $m$ dimensional process. 

Our aim throughout this paper is to estimate the $L_2$-distance $D^2:=\frac{1}{4\pi} \int_{-\pi}^\pi (f_{11}(\lambda)-f_{22}(\lambda))^2 d\lambda$, where $f_{11}(\lambda)$ and $f_{22}(\lambda)$ are the spectral densities of the first and the second process respectively. Under the null hypothesis
\begin{align}
\label{null}
\text{H}_0:\quad f_{11}(\lambda)=f_{22}(\lambda)
\end{align}
the distance $D^2$ equals zero while it is strictly positive if $f_{11}(\lambda) \not= f_{22}(\lambda)$ for $\lambda \in A$, where $A$ is a subset of $[-\pi,\pi]$ with positive Lebesgue measure. We will estimate $D^2$ by sums of the (squared) periodogram, where the sum goes over the Fourier coefficents of the smaller time series. Asymptotic normality both under the null and the alternative will be derived and since the variance terms can be easily estimated also under the alternative, asymptotic confidence intervalls and a precise hypothesis test can be constructed next to the test for \eqref{null} [see Remark \ref{rem1}]. Furthermore our approach has much wider application like testing for zero correlation, discriminant analysis or clustering of time series with unequal length [see Remark \ref{rem2}--\ref{remunkor}], and a simulation study will indicate that some of our assumptions are in fact not necessary (for example our method seems to work also for Long Memory processes). 

\section{The test statistic}
\def\theequation{2.\arabic{equation}}
\setcounter{equation}{0}
Let $n_1, n_2 \in \N$ with $n_1 \leq n_2$ and consider the two stationary time series
\begin{align}
\label{erst}
X_t^{(1)}&=\sum_{l=-\infty}^\infty \psi_l^{(1)} Z_{t-l}^{(1)} \quad t=1,...,n_1  \quad &X_t^{(2)}=\sum_{l=-\infty}^\infty \psi_l^{(2)} Z_{t-l}^{(2)} \quad t=1,...,n_2
\end{align}

where the $Z_t^{(j)}$ are independent and identically standard normal distributed for $j=1,2$ and
\begin{align}
\label{corr}
\E(Z_{t_1}^{(1)}Z_{t_2}^{(2)})= \begin{cases} \rho & \text{if } t_2= \lfloor t_1 q_{n_1,n_2} \rfloor-\lfloor q_{n_1,n_2}-1\rfloor \\ 0 & \text{else} \end{cases}  
\end{align}

where $q_{n_1,n_2}=\frac{n_2}{n_1}$ and $\rho \in [0,1]$. This roughly speaking means that changes in the time series with less observations influence the more frequently observed series but not vice versa, which is for example the case if interest rates and stock returns are compared. Throughout the paper we also assume that the technical condition
\begin{align}
\label{sumbed}
\sum_{l=-\infty}^\infty \psi_l^{(j)} |l|^\alpha <\infty 
\end{align}
is satisfied for an $\alpha > 1/2$ ($j=1,2$). Note that the assumption of Gaussianity is only imposed to simplify technical arguments [see Remark \ref{rem3}]. Furthermore innovations with variances different to $1$ can be included by choosing other coefficents $\psi_l^{(j)}$. We define the spectral densities $f_{ij}(\lambda)$ ($i,j=1,2$) through
\begin{eqnarray*}
f_{jj}(\lambda):= \frac{1}{2\pi} \left|\sum_{l=-\infty}^\infty \psi_l^{(j)} \exp(-i \lambda l)\right|^2, \quad \quad
f_{12}(\lambda):=\frac{\rho}{2\pi} \sum_{l,m=-\infty}^\infty \psi_l^{(1)}\psi_m^{(2)} \exp(-i \lambda (l-m)) = \overline{f_{21}(\lambda)}.
\end{eqnarray*}

An unbiased (but not consistent) estimator for $f_{jj}(\lambda)$ is given by the periodogram
\begin{align}
I_j(\lambda):=\frac{1}{2\pi n_j} \Bigl | \sum_{t=1}^{n_j} X_t^{(j)}\exp(-i \lambda t) \Bigr |^2
\end{align}

and although the periodogram does not estimate the spectral density consistently, a Riemann-sum over the Fourier coefficents of an exponentiated periodogram is (up to a constant) a consistent estimator for the corresponding integral over the exponentiated spectral density. For example, Theorem 2.1 in \cite{detkinsvet2010} yields
\begin{align}
\label{D1}
\hat D_{1,n_1}:=\frac{1}{n_1} \sum_{k=1}^{\lfloor \frac{n_1}{2} \rfloor } I_1^2(\lambda_{1,k}) \xrightarrow{\text{ }P \text{ }} \frac{1}{2\pi}\int_{-\pi}^\pi f_{11}^2(\lambda) d\lambda =: D_1
\end{align}

where $\lambda_{1,k}:= \frac{2\pi k}{n_1}$ ($k=1,...,\lfloor \frac{n_1}{2} \rfloor$) are the Fourier coefficents of the smaller time series $X_t^{(1)}$. If we can show that
\begin{align}
\label{D2}
\hat D_{2,n_1}&:=\frac{1}{n_1} \sum_{k=1}^{\lfloor \frac{n_1}{2} \rfloor } I_2^2(\lambda_{1,k}) \xrightarrow{\text{ }P \text{ }} \frac{1}{2\pi}\int_{-\pi}^\pi f_{22}^2(\lambda) d\lambda =: D_2, \\
\label{D12}
\hat D_{12,n_1}&:=\frac{1}{n_1}  \sum_{k=1}^{\lfloor \frac{n_1}{2} \rfloor-1 } I_1(\lambda_{1,k})I_2(\lambda_{1,k+1}) \xrightarrow{\text{ }P \text{ }} \frac{1}{4\pi}\int_{-\pi}^\pi f_{11}(\lambda) f_{22}(\lambda) d\lambda =: D_{12} ,
\end{align}

we can construct an consistent estimator for $D^2$ through $\hat D_{n_1}^2:=\frac{1}{2}(\hat D_{1,n_1}+\hat D_{2,n_1})-2\hat D_{12,n_1}$.
Although \eqref{D2} looks very much like \eqref{D1}, note that the convergence in \eqref{D2} is different since the coefficents $\lambda_{1,k}$ are not necessarily the Fourier coefficents of the time series $X_{t}^{(2)}$. This implies that the proof of \eqref{D2} has to be done in a completely different way than the proof of \eqref{D1} in \cite{detkinsvet2010}. We now obtain the following main theorem.

\begin{satz}
\label{thm1} If $f_{11}(\lambda)$, $f_{22}(\lambda)$ and $f_{12}(\lambda)$ are H\"older continuous of order $L>1/2$ and
\begin{align}
\label{bed}
\frac{n_2}{n_1} \rightarrow  Q
\end{align}
for a $Q \in \R$, then as $n_1 \rightarrow \infty$
\begin{eqnarray*}
\sqrt{n_1} \left( \hat D_{1,n_1}-D_1 , \hat D_{12,n_1}-D_{12}, \hat D_{2,n_1}-D_{2} \right)^T \xrightarrow{\text{ }D \text{ }} N(0,\left(\Sigma_{ij}\right)_{i,j=1}^3)
\end{eqnarray*}
with
\begin{eqnarray*}
\Sigma_{11} &=& 5 \int_{-\pi}^\pi f_{11}^4(\lambda) d\lambda, \quad  \Sigma_{12} = \int_{-\pi}^\pi f_{11}^3(\lambda)f_{22}(\lambda) d\lambda+\int_{-\pi}^\pi f_{11}^2(\lambda)|f_{12}(\lambda)|^2,
\\ \Sigma_{13} &=& \int_{-\pi}^\pi f_{12}^2(\lambda)f_{21}^2(\lambda)+4\int_{-\pi}^\pi f_{11}(\lambda)|f_{12}(\lambda)|^2f_{22}(\lambda),\quad \Sigma_{22} =  \frac{3}{4} \int_{-\pi}^\pi f_{11}^2(\lambda)f_{22}^2(\lambda) d\lambda+\frac{1}{2}\int_{-\pi}^\pi f_{11}(\lambda)|f_{12}(\lambda)|^2f_{22}(\lambda) d\lambda
\\ \Sigma_{23} &=& \int_{-\pi}^\pi f_{11}(\lambda)f_{22}^3(\lambda) d\lambda+\int_{-\pi}^\pi f_{22}^2(\lambda)|f_{12}(\lambda)|^2,\quad \Sigma_{33} = 5 \int_{-\pi}^\pi f_{22}^4(\lambda) d\lambda .
\end{eqnarray*}
\end{satz}

Although condition \eqref{bed} imposes some restrictions on the growth rate of $n_1$ and $n_2$, it is not very restrictive, since in practice there usually occur situations where even $n_2=Qn_1$ holds for a $Q \in \N$ (if for example daily data are compared with monthly data) and on the other hand this condition needs only to be satisfied in the limit. From Theorem $\ref{thm1}$ it now follows by a straightforward application of the Delta-Method that $\sqrt{n_1} (\hat D_{n_1}^2-D^2) \xrightarrow{\text{ }D \text{ }} N(0,\sigma^2)$, where $\sigma^2:= \frac{1}{\pi} \left \{ \frac{\Sigma_{11}+\Sigma_{33}}{4}+4\Sigma_{22}+\frac{\Sigma_{13}}{2}-2\Sigma_{12}-2\Sigma_{23} \right \}$, which becomes $\sigma_{\text{H}_0}^2= \frac{3}{2 \pi} \int_{-\pi}^{\pi} f_{11}^4(\lambda)d\lambda+\frac{1}{2\pi} \int_{-\pi}^{\pi} |f_{12}|^4 d\lambda$ under $\text{H}_0$. To obtain a consistent estimator for the variance under the null hypothesis we define
\begin{align*}
I_{12}(\lambda):&= \frac{1}{2\pi \sqrt{n_1n_2}} \sum_{p_1=1}^{n_1} X_{p_1}^{(1)} \exp(-i \lambda p_1) \sum_{p_2=1}^{n_2} X_{p_2}^{(2)} \exp(i \lambda p_2)= \overline{I_{21}(\lambda)}
\end{align*}
and analogous to the proof of Theorem \ref{thm1} it can be shown that
\begin{align*}
\hat \sigma_{\text{H}_0}^2 :=\frac{1}{4 n_1} \sum_{k=1}^{\lfloor \frac{n_1}{2} \rfloor } \Bigl(I_1^4(\lambda_{1,k})+I_2^4(\lambda_{1,k}) \Bigr) + Re \Bigl( \frac{1}{2n_1} \sum_{k=1}^{\lfloor n_1/2 \rfloor -1} I_{12}^2(\lambda_{1,k})I_{21}^2(\lambda_{1,k+1}) \Bigr) \xrightarrow{\text{ }P \text{ }} \sigma_{\text{H}_0}^2 .
\end{align*}

Therefore an asymptotic niveau-$\alpha$-test for \eqref{null} is given by: reject \eqref{null} if
\begin{align}
\label{test}
\sqrt{n_1} \frac{\hat D_{n_1}^2}{\sqrt{\hat \sigma_{\text{H}_0}^2}} > u_{1-\alpha} ,
\end{align}

where $u_{1-\alpha}$ denotes the $(1-\alpha)$ quantile of the standard normal distribution. This test has asymptotic power \linebreak $\Phi \left(\sqrt{n_1} \frac{D^2}{\sigma}- \frac{\sqrt{\hat \sigma_{\text{H}_0}^2}}{\sigma} u_{1-\alpha} \right)$ where $\Phi$ is the distribution function of the standard normal distribution. This yields that the test \eqref{test} has asymptotic power one for all alternatives with $D^2 >0$. 

\begin{rem}
\label{rem1} ~~~~ \\
{\rm  It is straightforward to construct an estimator $\hat{\sigma}^2$, which converges to the variance $\sigma^2$ also under the alternative. This enables us to construct asymptotic $(1-\alpha)$ confidence intervals for $D^2$. The same statement holds, if we consider the normalized measure $R^2:= \frac{2D^2}{D_1+D_2}$,which can be estimated by $\hat R_{n_1}^2:= \frac{2 \hat D_{n_1}^2}{ \hat D_{1,n_1} + \hat D_{2,n_1}}$.

From Theorem \ref{thm1} and a straightforward application of the delta method, it follows that
\begin{align}
\label{zentr2}
\sqrt{n_1} (\hat R_{n_1}^2-R^2) \xrightarrow{\text{ }D \text{ }} N(0,\sigma_1^2),
\end{align}
where $\sigma_1^2$ can be easily calculated. By considering a consistent estimator $\hat{\sigma}_1^2$ for $\sigma_1^2$ (which can be constructed through the corresponding Riemann sums of the periodogram), \eqref{zentr2} provides an asymptotic level $\alpha$ test for the so called \textit{precise hypothesis}
\begin{equation} \label{precise}
H_0 : R^2  > \varepsilon \quad \mbox {versus} \quad  H_1 :
R^2 \leq \varepsilon \:,
\end{equation}
where $\varepsilon >0$ [see \cite{bergdela1987}]. This hypothesis is of interest, because spectral densities of time series in real-world applications are usually never exactly equal and a more realistic question is then to ask, if the processes have approximately the same spectral measure. An asymptotic level $\alpha$ test for \eqref{precise} is obtained by rejecting the null hypothesis, whenever $\hat R_{n_1}^2 - \varepsilon < \frac{ \hat \sigma_1 }{\sqrt{n_1}} u_{\alpha}$.

}
\end{rem}

\begin{rem}
\label{rem2} ~~~~ \\
{\rm
Theorem \ref{thm1} can also be employed for a cluster and a discriminant analysis of time series data with different length, since it yields an estimator for the distance measure $d(f_{11},f_{22})$, where
\begin{eqnarray*}
d(f,g) &=& \left( 1- \frac{2\int_{-\pi}^\pi f(\lambda)g(\lambda) d\lambda}{\int_{-\pi}^\pi f^2(\lambda)d\lambda + \int_{-\pi}^\pi g^2(\lambda) d\lambda} \right)^{1/2},
\end{eqnarray*}
which can take values between 0 and 1. A value close to 0 indicates some kind of similarities between two processes, whereas a value close to 1 exhibits dissimilarities in the second-order structure. The distance measure $d(f_{11},f_{22})$ can be estimated by
\begin{eqnarray}
\hat{d}_{12} &=& \left( \max \left(1- \frac{2 \hat{D}_{12,n_1}}{\hat{D}_{1,n_1} + \hat{D}_{2,n_1}},0 \right) \right)^{1/2} , \label{distance}
\end{eqnarray}
where the maximum is necessary, because the term $1- \frac{2 \hat{D}_{12,n_1}}{\hat{D}_{1,n_1} + \hat{D}_{2,n_1}}$ might be negative.
}
\end{rem}

\begin{rem}
\label{remunkor} ~~~~ \\
{\rm The main ideas of the proof of Theorem \ref{thm1} can be furthermore employed to construct tests for various other hypothesis. For example a test for zero correlation can be derived by testing for $ f_{12} \equiv 0$ which can be done by estimating $\int_{-\pi}^\pi |f_{12}(\lambda)|^2 d\lambda$. An estimator for this quantity is easily derived using the above approach and furthermore the calculation of the variance is straightforward, which we omit for the sake of brevity.
}
\end{rem}

\begin{rem}
\label{rem3} ~~~~ \\
{\rm Although we only considered the bivariate case, our method can be easily extended to an $m$ dimensional process. Moreover, a cumbersome but straightforward examination yields that our test also has asymptotic level $\alpha$, if we skip the assumption of Gaussianity since (under the null hypothesis) all terms which consist the fourth cumulants of the processes $Z^{(i)}$ cancel out. A similar phenomenon can be observed for the tests proposed by \cite{eichler2008}, \cite{detkinsvet2010}, \cite{dethil2011} and \cite{detprevet2010}.
}

\end{rem}

\section{Finite sample study}
\def\theequation{3.\arabic{equation}}
\setcounter{equation}{0}
\subsection{Size and power of the test}
In this section we study the size and the power of test \eqref{test} in the case of finite samples. All simulations are based on 1000 iterations and we consider all different combinations of $n_1,n_2 \in \{256,384,512,640\}$ with $n_1 \le n_2$. For the sake of brevity we only present the results for the case $\rho = 0$ and note that the rejection frequencies do not change at all if we consider correlations different to zero. We furthermore tested our approach using non-linear GARCH models and obtained a very good performance also in this case. The results are not displayed for the sake of brevity but are available from the authors upon request. To demonstrate the approximation of the nominal level, we consider the five processes
\begin{eqnarray*}
\textbf{X}^1: X_t &=&Z_t, \quad \textbf{X}^2: X_t=-0.8X_{t-1}+Z_t, \quad \textbf{X}^3: X_t= Z_t-0.8Z_{t-1}, \quad \textbf{X}^4: X_t\sim FARIMA(0.45,0,0.8),
\\   \textbf{X}^5: X_t&=& Z_t 1_{[t\le0.5T]}+0.8X_{t-1}1_{[0.5T\le t\le 0.75T]}+Z_t 1_{[t\ge0.75T]} \quad \text{for } t=1,...,T,
\end{eqnarray*}
where the $FARIMA(0.45,0,0.8)$-model corresponds to a LongMemory-process given by  $(1-B)^{0.45} X_t=(1-0.8B) Z_t$ with the backshift-operator $B$ (i.e. $B^jX_t=X_{t-j}$). Note that the models $\textbf{X}^4$ and $\textbf{X}^5$ both do not fit into the theoretical framework considered in section 2, since for the $FARIMA(0.45,0,0.8)$-process we obtain $\sum_{l=-\infty}^\infty |\psi_l|=\infty $ which contradicts \eqref{sumbed} and the structural-break model $\textbf{X}^5$ does not even has a stationary solution. Nevertheless since these models are of great interest in practice, we investigate the performance of our approach in these cases as well. The results are given in Table \ref{tab1} and it can be seen that the test \eqref{test} is very robust against different choices of $n_1$ and $n_2$. Furthermore our method also seems to work for the models $\textbf{X}^4$ and $\textbf{X}^5$ although the convergence is slightly slower. 

To study the power of the test we additionally present the results of a comparison of $\textbf{X}^3$ with $\textbf{X}^j$ for $j \in \{1,2\}$ and $\textbf{X}^1$ with $\textbf{X}^5$ (all other comparison between the processes yield better results than the depicted ones).

\subsection{Real world data}
\begin{figure}[hbt]
	\centering
	\includegraphics[width=12cm]{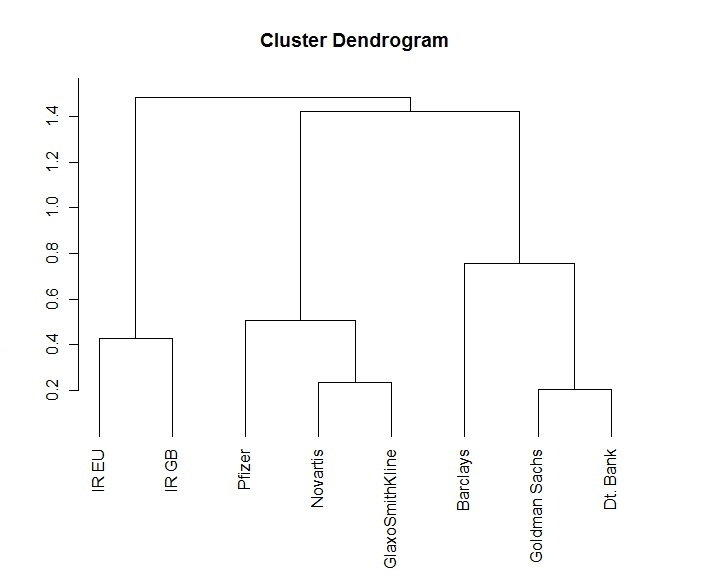}
	\caption{Clustering of financial time series data.} \label{cluster}
\end{figure}

In this section we investigate how the clustering-method described in Remark \ref{rem2} performs, if it is applied to real world data. Therefore we took three log-returns of stock prices from the financial sector, three  log-returns from the health sector and two key interest rates. Exemplarily for the finance sector we choosed the stocks of Barclays, Deutsche Bank and Goldman Sachs and the health sector is represented by GlaxoSmithKline, Novartis and Pfizer. The key interest rates were taken from Great Britain and the EU and all time series data were recorded between March 1st, 2003 and July 29th, 2011. While the interest rates data were observed monthly, the stock prices were recorded daily or weekly . However, even if two stock prices were observed daily they might differ in length, since they are for example traded on different stock exchanges with different trading days. The result of our cluster analysis using \eqref{distance} is presented in the dendrogram given in Figure \ref{cluster}. We get three different groups which correspond to the finance sector, the health sector and the key interes rates.

\bigskip

{\bf Acknowledgements}
This work has been supported in part by the
Collaborative Research Center ``Statistical modeling of nonlinear
dynamic processes'' (SFB 823, Teilprojekt C1, C4) of the German Research Foundation
(DFG).

\section{Appendix: Technical details}
\def\theequation{4.\arabic{equation}}
\setcounter{equation}{0}
{\bf Proof of Theorem \ref{thm1}:} By using the Cramer-Wold device, we have to show that
\begin{align*}
c^T \sqrt{n_1} \Bigl \{ (\hat D_{1,n_1}, \hat D_{12,n_1},\hat D_{2,n_1})^T-(D_1,D_{12}, D_{2})^T \Bigr \} \xrightarrow{\text{ }D \text{ }} N(0, c^T \Sigma c) 
\end{align*}

for all vectors $c \in \R^3$. For the sake of brevity, we restrict ourselve to the case $c=(0,1,0)^T$ since the more general follows with exactly the same arguments. Therefore we show $\hat T_{n_1}:=\sqrt{n_1} (\hat D_{12,n_1}-D_{12}) \xrightarrow{\text{ }D \text{ }} N(0,\Sigma_{22})$ and we do that by using the method of cumulants, which is described in chapter 2.3. of \cite{brillinger1981} (and whose notations we will make heavy use of), i.e. in the following it is proved that
\begin{align}
\label{cum1}
&\cum_l(\hat T_{n_1})=o(1) \quad \mbox{for $l=1$ and $l\geq 3$}, \\
\label{cum2}
&\cum_2(\hat T_{n_1}) \xrightarrow{n_1 \rightarrow \infty} \Sigma_{22}, 
\end{align}
which will yield the assertion. 

{\bf Proof of \eqref{cum1} for the case $l=1$:}  Because of the symmetry of the periodogram, it is 
\begin{eqnarray*}
\E(\hat D_{12,n_1}) &=& \frac{1}{2n_1} \sum_{k=-\lfloor \frac{n_1-1}{2} \rfloor}^{\lfloor \frac{n_1}{2} \rfloor } \frac{1}{(2\pi)^2 n_1n_2} \sum_{p_1,q_1=1}^{n_1} \sum_{p_2,q_2=1}^{n_2} \sum_{l_1,m_1=-\infty}^{\infty}\sum_{l_2,m_2=-\infty}^{\infty} \psi_{l_1}^{(1)}\psi_{m_1}^{(1)}\psi_{l_2}^{(2)}\psi_{m_2}^{(2)}
\\ && E(Z_{p_1-l_1}^{(1)} Z_{q_1-m_1}^{(1)} Z_{p_2-l_2}^{(2)} Z_{q_2-m_2}^{(2)}) e^{-i\lambda_{1,k}(p_1-q_1)-i\lambda_{1,k+1}(p_2-q_2)} + O(1/n_1)
\end{eqnarray*}

and because of the standard normality of the innovations we obtain
\begin{align*}
\E(Z_{p_1-l_1}^{(1)} Z_{q_1-m_1}^{(1)} Z_{p_2-l_2}^{(2)} Z_{q_2-m_2}^{(2)})=&\E(Z_{p_1-l_1}^{(1)} Z_{q_1-m_1}^{(1)})\E( Z_{p_2-l_2}^{(2)} Z_{q_2-m_2}^{(2)})+\E(Z_{p_1-l_1}^{(1)}Z_{q_2-m_2}^{(2)}) \E( Z_{q_1-m_1}^{(1)} Z_{p_2-l_2}^{(2)} ) \\
&+\E(Z_{p_1-l_1}^{(1)} Z_{p_2-l_2}^{(2)} ) \E(Z_{q_1-m_1}^{(1)} Z_{q_2-m_2}^{(2)})
\end{align*}

which yields that $\E(\hat D_{12,n_1})$ (without the $O(1/n_1)$-term) can be divided into the sums of three terms which are called $A$, $B$ and $C$ respectively. For the first term we obtain the conditions $p_1=q_1+l_1-m_1, p_2=q_2+l_2-m_2$ (all others cases are equal to zero). This results in
\begin{eqnarray}
\nonumber
A &=& \frac{1}{2n_1} \sum_{k=-\lfloor \frac{n_1-1}{2} \rfloor}^{\lfloor \frac{n_1}{2} \rfloor } \frac{1}{(2\pi)^2 n_1n_2} \sum_{l_1,l_2,m_1,m_2 = -\infty}^\infty \sum_{\substack {q_1=1 \\ 1\le q_1+l_1-m_1\le n_1} }^{n_1}\sum_{\substack {q_2=1 \\ 1\le q_2+l_2-m_2\le n_2} }^{n_2}\psi_{l_1}^{(1)}...\psi_{m_2}^{(2)}e^{-i\lambda_{1,k}(l_1-m_1)-i\lambda_{1,k+1}(l_2-m_2)}
\\ \nonumber &=& \frac{1}{2n_1} \sum_{k=-\lfloor \frac{n_1-1}{2} \rfloor}^{\lfloor \frac{n_1}{2} \rfloor } \frac{1}{(2\pi)^2 n_1n_2} \sum_{q_1=1}^{n_1} \sum_{q_2=1}^{n_2}\sum_{l_1,l_2,m_1,m_2 = -\infty}^\infty\psi_{l_1}^{(1)}...\psi_{m_2}^{(2)}e^{-i\lambda_{1,k}(l_1-m_1)-i\lambda_{1,k+1}(l_2-m_2)} + o\left( \frac{1}{\sqrt{n}} \right) ,
\end{eqnarray}

where the last equality follows from
\begin{align}
\label{warumbedweglassen}
\frac{1}{n_j}\sum_{l: |l| < M n_j} \psi_l^{(j)} |l| = \frac{1}{n_j} \sum_{l: |l| < M n_j} \psi_l^{(j)} |l|^\alpha |l|^{1-\alpha} =o (1/n_j^\alpha)
\end{align}

with $M \in \R$, where \eqref{sumbed} was used. It now follows by the H\"older continuity condition that $A$ equals $ \frac{1}{2\pi} \int_{-\pi}^\pi f_{11}(\lambda)f_{22}(\lambda) d\lambda + o\left( \frac{1}{\sqrt{n}} \right)$. If we consider the summand $B$, we obtain the conditions $q_1 = \lfloor (p_2-l_2) q_{n_1,n_2} \rfloor +m_1 -\lfloor q_{n_1,n_2} -1 \rfloor$, \\ 
$\quad q_2 = \lfloor (p_1-l_1) q_{n_1,n_2} \rfloor +m_2-\lfloor q_{n_1,n_2} -1 \rfloor$ which yields
\begin{eqnarray*}
B &=& \frac{\rho^2}{2n_1} \sum_{k=-\lfloor \frac{n_1-1}{2} \rfloor}^{\lfloor \frac{n_1}{2} \rfloor } \frac{1}{(2\pi)^2 n_1n_2} \sum_{l_1,m_1, l_2,m_2=-\infty}^{\infty} \sum_{\substack{p_1=1 \\ 1 \leq \lfloor (p_1-l_1) q_{n_1,n_2} \rfloor +m_2-\lfloor q_{n_1,n_2} -1 \rfloor \leq n_2}}^{n_1} \sum_{\substack{p_2=1 \\ 1 \leq \lfloor (p_2-l_2) q_{n_1,n_2} \rfloor +m_1 -\lfloor q_{n_1,n_2} -1 \rfloor \leq n_1}}^{n_2} \\ && 
  \psi_{l_1}^{(1)}...\psi_{m_2}^{(2)}
 e^{-i\lambda_{1,k}(p_1-\lfloor (p_2-l_2) q_{n_1,n_2} \rfloor +m_1 -\lfloor q_{n_1,n_2} -1 \rfloor) } e^{-i\lambda_{1,k+1}(p_2-\lfloor (p_1-l_1) q_{n_1,n_2} \rfloor +m_2-\lfloor q_{n_1,n_2} -1 \rfloor ) } .
\end{eqnarray*}

If we now employ the identity
\begin{align}
\label{help1}
\frac{1}{n_1} \sum_{k=-\lfloor \frac{n_1-1}{2} \rfloor}^{\lfloor \frac{n_1}{2} \rfloor } e^{-i\lambda_{1,k} p}= \begin{cases} 1 & \text{ if } p=0, \pm n_1, \pm 2n_1, ... \\ 0 & \text{ else} \end{cases} ,
\end{align}
it follows with \eqref{bed} that if $p_1$ is chosen there are only finitely many $p_2$ which yields a non-zero summand. Therefore we obtain that
$B=o(1/\sqrt{n_1})$ and with the same arguments it can be shown that $C=o(1/\sqrt{n_1})$. $\hfill \Box$

{\bf Proof of \eqref{cum2}:} It is
\begin{align*}
\cum_2(\sqrt{n_1} \hat{D}_{12,n_1})&=\frac{1}{n_1}  \sum_{k=1}^{\lfloor \frac{n_1}{2} \rfloor -1} \cum_2(I_1(\lambda_{1,k})I_2(\lambda_{1,k+1})) + \frac{1}{n_1}  \sum_{\substack{k_1,k_2=1 \\ k_1 \not= k_2}}^{\lfloor \frac{n_1}{2} \rfloor -1} \cum(I_1(\lambda_{1,k_1})I_2(\lambda_{1,k_1+1}),I_1(\lambda_{1,k_2})I_2(\lambda_{1,k_2+1}))
\end{align*}

and the assertion follows if we show that
\begin{eqnarray}
\label{first}  \frac{1}{n_1}  \sum_{k=1}^{\lfloor \frac{n_1}{2} \rfloor -1} \cum_2(I_1(\lambda_{1,k})I_2(\lambda_{1,k+1}))  &\xrightarrow{n_1 \rightarrow \infty}& \frac{3}{4\pi} \int_{-\pi}^\pi f_{11}^2(\lambda)f_{22}^2(\lambda) d\lambda, \\
\label{second}  \frac{1}{n_1}  \sum_{\substack{k_1,k_2=1 \\ k_1 \not= k_2}}^{\lfloor \frac{n_1}{2} \rfloor -1} \cum(I_1(\lambda_{1,k_1})I_2(\lambda_{1,k_1+1}),I_1(\lambda_{1,k_2})I_2(\lambda_{1,k_2+1}))  &\xrightarrow{n_1 \rightarrow \infty}& \frac{1}{2\pi} \int_{-\pi}^\pi  f_{11}(\lambda)|f_{12}(\lambda)|^2f_{22}(\lambda) d\lambda .
\end{eqnarray}

We present a detailed proof of \eqref{first} and then comment briefly on \eqref{second} since it is proved analogously. Employing the symmetry of the periodogram again, we get
\begin{eqnarray}
 &&\frac{1}{n_1} \sum_{k=1}^{\lfloor \frac{n_1}{2} \rfloor -1}  \cum_2(I_1(\lambda_{1,k})I_2(\lambda_{1,k+1}))=  \frac{1}{2n_1} \sum_{k=-\lfloor \frac{n_1-1}{2} \rfloor}^{\lfloor \frac{n_1}{2} \rfloor } \cum_2(I_1(\lambda_{1,k})I_2(\lambda_{1,k+1})) +O(1/n_1) \nonumber \\ \nonumber
 &=&  \frac{1}{2n_1} \sum_{k=-\lfloor \frac{n_1-1}{2} \rfloor}^{\lfloor \frac{n_1}{2} \rfloor }  \frac{1}{(2\pi)^4n_1^2n_2^2} \sum_{j=1}^2  \sum_{p_j,q_j,r_j,s_j=1}^{n_j} \sum_{a_j,b_j,c_j,d_j=-\infty}^\infty \psi_{a_1}^{(1)}...\psi_{d_2}^{(2)}
\\ \nonumber && e^{-i\lambda_{1,k}(p_1-q_1+r_1-s_1)-i\lambda_{1,k+1}(p_2-q_2+r_2-s_2)} \cum(Z_{p_1-a_1}^{(1)} Z_{q_1-b_1}^{(1)} Z_{p_2-a_2}^{(2)} Z_{q_2-b_2}^{(2)}, Z_{r_1-c_1}^{(1)} Z_{s_1-d_1}^{(1)} Z_{r_2-c_2}^{(2)} Z_{s_2-d_2}^{(2)}) +O(1/n_1) \\ \label{varianz}
&=& \sum_\nu  \frac{1}{2n_1} \sum_{k=-\lfloor \frac{n_1-1}{2} \rfloor}^{\lfloor \frac{n_1}{2} \rfloor }  \frac{1}{(2\pi)^4n_1^2n_2^2} \sum_{j=1}^2  \sum_{p_j,q_j,r_j,s_j=1}^{n_j} \sum_{a_j,b_j,c_j,d_j=-\infty}^\infty \psi_{a_1}^{(1)}...\psi_{d_2}^{(2)}
\\ \nonumber && e^{-i\lambda_{1,k}(p_1-q_1+r_1-s_1)-i\lambda_{1,k+1}(p_2-q_2+r_2-s_2)}  \cum(Z_i^{(j)}; (i,j) \in \nu_1) \cdots \cum(Z_i^{(j)}; (i,j) \in \nu_4) +O(1/n_1)
\end{eqnarray}

where the sum goes over all indecomposable partitions $\nu=\nu_1 \cup ... \cup \nu_4$ of
\begin{align*}
\begin{matrix}
Z_{p_1-a_1}^{(1)} &  Z_{q_1-b_1}^{(1)} & Z_{p_2-a_2}^{(1)} & Z_{q_2-b_2}^{(1)}  \\
Z_{r_1-c_1}^{(2)} &  Z_{s_1-d_1}^{(2)} & Z_{r_2-c_2}^{(2)} & Z_{s_2-d_2}^{(2)}  
\end{matrix}
\end{align*}

with $|\nu_i|=2$ $\forall i=1,...,4$ (we only have to consider partitions with two elements in each set, because of the Gaussianity of the innovations; in the non-Gaussian case we would get an additional term containing the fourth cumulant). Every chosen partition will imply conditions for the choice of $p_j,q_j,r_j,s_j$ as in the calculation of the expectation. For some partitions there will not be a $p_j,q_j,r_j,s_j$ in the exponent of $e$ after inserting the conditions and for other partitions there will still remain one. Let us take an example of the latter one and consider the partition which corresponds to $\cum(Z_{p_1-a_1}^{(1)},Z_{s_1-d_1}^{(1)}) \cum(Z_{q_1-b_1}^{(1)},Z_{r_1-c_1}^{(1)}) \cum(Z_{p_2-a_2}^{(2)},Z_{r_2-c_2}^{(2)}) \cum(Z_{q_2-b_2}^{(2)},Z_{s_2-d_2}^{(2)})$. We name the corresponding term of this partition in \eqref{varianz} with $V_2$ and obtain the conditions $p_1 =s_1+a_1-d_1$,  $q_1=r_1+b_1-c_1$, $p_2=r_2+a_2-c_2$, $q_2=s_2+b_2-d_2$ which yields
\begin{eqnarray*}
V_2 &=& \frac{1}{2n_1} \sum_{k=-\lfloor \frac{n_1-1}{2} \rfloor}^{\lfloor \frac{n_1}{2} \rfloor } \frac{1}{(2\pi)^4n_1^2n_2^2} \sum_{j=1}^2 \sum_{a_j,b_j,c_j,d_j=-\infty}^\infty\sum_{\substack{s_1,r_1=1\\ 1 \le s_1+a_1-d_1 \le n_1 \\ 1 \le r_1+b_1-c_1\le n_1}}^{n_1}\sum_{\substack{r_2,s_2=1\\ 1 \le r_2+a_2-c_2 \le n_2 \\ 1 \le s_2+b_2-d_2 \le n_2}}^{n_2} \psi_{a_1}^{(1)}...\psi_{d_2}^{(2)}
\\ && e^{-i\lambda_{1,k}(a_1-d_1+c_1-b_1)-i\lambda_{1,k+1}(2r_2-2s_2+a_2-c_2+d_2-b_2)} \\
&=& \frac{1}{2n_1} \sum_{k=-\lfloor \frac{n_1-1}{2} \rfloor}^{\lfloor \frac{n_1}{2} \rfloor }  \frac{1}{(2\pi)^4n_1^2n_2^2} \sum_{j=1}^2 \sum_{a_j,b_j,c_j,d_j=-\infty}^\infty\sum_{s_1,r_1=1}^{n_1}\sum_{r_2,s_2=1}^{n_2} \psi_{a_1}^{(1)}...\psi_{d_2}^{(2)}
\\ && e^{-i \lambda_{1,k+1}(2r_2-2s_2)}e^{-i\lambda_{1,k}(a_1-d_1+c_1-b_1)-i\lambda_{1,k+1}(a_2-c_2+d_2-b_2)}+o(1/\sqrt{n_1}) ,
\end{eqnarray*}

where the last equality again follows with \eqref{warumbedweglassen}. Now as in the handling of $B$ in the calculation of the expectation, \eqref{help1} implies that $V_2=o(1)$. \\
Every other indecomposable partition is treated in exactly the same way and there are only three partitions which corresponding term in \eqref{varianz} does not vanish in the limit. These partitions correspond to one of the following three terms:
\begin{align*}
&1) \quad \cum(Z_{p_1-a_1}^{(1)} , Z_{q_1-b_1}^{(1)}) \cum(Z_{r_1-c_1}^{(1)} , Z_{s_1-d_1}^{(1)} ) \cum(Z_{p_2-a_2}^{(2)} , Z_{s_2-d_2}^{(2)}  ) \cum(Z_{q_2-b_2}^{(2)}  , Z_{r_2-c_2}^{(2)}  ) \\
&2) \quad \cum(Z_{p_1-a_1}^{(1)} , Z_{s_1-d_1}^{(1)}) \cum(Z_{q_1-b_1}^{(1)} , Z_{r_1-c_1}^{(1)} ) \cum(Z_{p_2-a_2}^{(2)} , Z_{q_2-b_2}^{(2)}  ) \cum( Z_{r_2-c_2}^{(2)}  ,  Z_{s_2-d_2}^{(2)} ) \\
&3) \quad \cum(Z_{p_1-a_1}^{(1)} , Z_{s_1-d_1}^{(1)}) \cum(Z_{q_1-b_1}^{(1)} , Z_{r_1-c_1}^{(1)} ) \cum(Z_{p_2-a_2}^{(2)} , Z_{s_2-d_2}^{(2)}  ) \cum(Z_{q_2-b_2}^{(2)}  , Z_{r_2-c_2}^{(2)}  )
\end{align*}

We will exemplarily present the calculation concerning the $1)$ partition and denote the corresponding sum in \eqref{varianz} with $V_1$. We get that $V_1$ equals
\begin{eqnarray*}
&&  \frac{1}{2n_1} \sum_{k=-\lfloor \frac{n_1-1}{2} \rfloor}^{\lfloor \frac{n_1}{2} \rfloor } \frac{1}{(2\pi)^4n_1^2n_2^2} \sum_{j=1}^2 \sum_{a_j,b_j,c_j,d_j=-\infty}^\infty\sum_{\substack{q_1,s_1=1\\ 1 \le q_1+a_1-b_1 \le n_1 \\ 1 \le s_1+c_1-d_1\le n_1}}^{n_1}\sum_{\substack{r_2,s_2=1\\ 1 \le s_2+a_2-d_2 \le n_2 \\ 1 \le r_2+b_2-c_2 \le n_2}}^{n_2} \psi_{a_1}^{(1)}...\psi_{d_2}^{(2)}e^{-i\lambda_{1,k}(a_1-b_1+c_1-d_1)-i\lambda_{1,k+1}(a_2-d_2+c_2-b_2)}  \\
&=&  \frac{1}{2n_1} \sum_{k=-\lfloor \frac{n_1-1}{2} \rfloor}^{\lfloor \frac{n_1}{2} \rfloor }\frac{1}{(2\pi)^4n_1^2n_2^2} \sum_{j=1}^2 \sum_{a_j,b_j,c_j,d_j=-\infty}^\infty\sum_{q_1,s_1=1}^{n_1}\sum_{r_2,s_2=1}^{n_2} \psi_{a_1}^{(1)}...\psi_{d_2}^{(2)} e^{-i\lambda_{1,k}(a_1-b_1+c_1-d_1)-i\lambda_{1,k+1}(a_2-d_2+c_2-b_2)} +o(1/\sqrt{n_1})
\end{eqnarray*}

by using \eqref{warumbedweglassen}. Now the H\"older continuity condition implies $V_1 = \frac{1}{4 \pi} \int_{-\pi}^\pi f_{11}(\lambda)^2 f_{22}(\lambda)^2 d\lambda +o(1/\sqrt{n_1})$ and since the partitions $2)$ and $3)$ yield the same result, we have shown \eqref{first}. \\
With the same arguments as in the proof of \eqref{first} it can be seen that
\begin{align*}
&\frac{1}{n_1}  \sum_{\substack{k_1,k_2=1 \\ k_1 \not= k_2}}^{\lfloor \frac{n_1}{2} \rfloor -1} \cum(I_1(\lambda_{1,k_1})I_2(\lambda_{1,k_1+1}),I_1(\lambda_{1,k_2})I_2(\lambda_{1,k_2+1})) \\
=& \frac{1}{n_1}  \sum_{k_1=1}^{\lfloor \frac{n_1}{2} \rfloor -1} \cum(I_1(\lambda_{1,k_1})I_2(\lambda_{1,k_1+1}),I_1(\lambda_{1,k_1+1})I_2(\lambda_{1,k_1+2})) \\
&+\frac{1}{n_1}  \sum_{k_1=1}^{\lfloor \frac{n_1}{2} \rfloor -1} \cum(I_1(\lambda_{1,k_1})I_2(\lambda_{1,k_1+1}),I_1(\lambda_{1,k_1-1})I_2(\lambda_{1,k_1})) +o(1)
\end{align*}

and it is shown completely analogously to the proof of \eqref{first} that 
\begin{align*}
\frac{1}{n_1}  \sum_{k_1=1}^{\lfloor \frac{n_1}{2} \rfloor -1} \cum(I_1(\lambda_{1,k_1})I_2(\lambda_{1,k_1+1}),I_1(\lambda_{1,k_1+1})I_2(\lambda_{1,k_1+2})), \quad
\frac{1}{n_1}  \sum_{k_1=1}^{\lfloor \frac{n_1}{2} \rfloor -1} \cum(I_1(\lambda_{1,k_1})I_2(\lambda_{1,k_1+1}),I_1(\lambda_{1,k_1-1})I_2(\lambda_{1,k_1}))
\end{align*}
both converge to $\frac{1}{4\pi} \int_{-\pi}^\pi  f_{11}(\lambda)|f_{12}(\lambda)|^2f_{22}(\lambda) d\lambda$, which yields \eqref{second}.  $\hfill \Box$

{\bf Proof of \eqref{cum1} for the case $l \geq 3$:} Since the proof is done by combining standard cumulants methods with the arguments that are used in the previous proof, we will restrict ourselve to a brief explanation of the main ideas. We obtain
\begin{align*}
\cum_l(\sqrt{n_1} D_{12,n_1})=&\frac{1}{(2n_1)^{l/2}} \sum_{j_1=1}^l \sum_{k_{j_1}=-\lfloor \frac{n_1-1}{2} \rfloor}^{\lfloor \frac{n_1}{2} \rfloor} \frac{1}{(2\pi)^{2l}n_1^ln_2^l} \sum_{j_2=1}^2 \sum_{a_{j_1,j_2},b_{j_1,j_2}=-\infty}^\infty \sum_{p_{j_1,j_2},q_{j_1,j_2}=1}^{n_{j_2}} \psi_{a_{1,1}}^{(1)} \cdots \psi_{b_{l,2}}^{(2)} \\
& \exp(-i \lambda_{1,k}(p_{11}-q_{11})-i \lambda_{1,k+1}(p_{12}-q_{12}))\cdots \exp(-i \lambda_{1,k}(p_{l1}-q_{l1})-i \lambda_{1,k+1}(p_{l2}-q_{l2})) \\
& \cum(Z_{p_{11}-a_{11}}^{(1)}Z_{q_{11}-b_{11}}^{(1)} Z_{p_{12}-a_{12}}^{(2)}Z_{q_{12}-b_{12}}^{(2)},...,Z_{p_{l1}-a_{l1}}^{(1)}Z_{q_{l1}-b_{l1}}^{(1)} Z_{p_{l2}-a_{l2}}^{(2)}Z_{q_{l2}-b_{l2}}^{(2)})
\end{align*}

and if we now take a indecomposable partition of 
\begin{align*}
\begin{matrix}
Z_{p_{11}-a_{11}}^{(1)} & Z_{q_{11}-b_{11}}^{(1)} & Z_{p_{12}-a_{12}}^{(2)} & Z_{q_{12}-b_{12}}^{(2)} \\
\vdots & \vdots & \vdots & \vdots \\
Z_{p_{l1}-a_{l1}}^{(1)} & Z_{q_{l1}-b_{l1}}^{(1)} & Z_{p_{l2}-a_{l2}}^{(2)} & Z_{q_{l2}-b_{l2}}^{(2)}
\end{matrix}
\end{align*}

which consists only of sets with two elements (again this suffices because of the Gaussianity of the innovations), it follows directly that at most $2l$ of the $4l$ variables $p_{j_1,j_2}, q_{j_1,j_2}$ ($j_1=1,...,l$, $j_2=1,2$) are free to choose. By using the same arguments as in the calculation of the variance and the expectation it then follows by the indecomposability of the partition that in fact only $l+1$ of the remaining $2l$ variables $p_{j_1,j_2}, q_{j_1,j_2}$ are free to choose. This implies $\cum_l(\sqrt{n_1} D_{12,n_1})=O(n_1^{1-l/2})$ which yields the assertion. $\hfill \Box$

\bibliographystyle{apalike}

\medskip

\renewcommand{\baselinestretch}{1.15}

\begin{table}[htb!]
\small
\centering
\begin{tabular}{|| c|c|c|c  c  c|| c   c || c c || c  ||}
\hline
$n_1$ &$n_2$ & $\alpha$&$\textbf{X}^1$ & $\textbf{X}^2$     & $\textbf{X}^3$ &$\textbf{X}^1\textbf{X}^3$&$\textbf{X}^2\textbf{X}^3$&$\textbf{X}^4$ & $\textbf{X}^5$ &$\textbf{X}^1\textbf{X}^5$\\ \hline
 256&256 &0.05&0.041 &0.039& 0.045&0.773&0.628&0.053 & 0.053&0.523 \\
 & &0.1&0.088&0.128 &0.106&0.894 &0.875&0.117&0.131&0.685\\
 \hline
 256& 384 &0.05&0.049&0.049& 0.034&0.758 &0.619&0.047 & 0.054&0.551\\
 & &0.1&0.106&0.128 &0.099&0.877&0.841&0.114&0.147&0.719\\
 \hline
 256&512 &0.05&0.043&0.047&0.026&0.776&0.650& 0.045&0.031&0.539 \\
 & &0.1&0.085&0.139 &0.069&0.892&0.848&0.109&0.126&0.739\\
\hline
 256&640 &0.05&0.050&0.044&0.030&0.777&0.636& 0.037& 0.046 &0.563\\
 & &0.1&0.109&0.112& 0.081&0.904&0.859&0.097&0.122&0.755\\
 \hline
 384& 384 &0.05&0.036&0.047&0.037 &0.920&0.804& 0.054& 0.043&0.693\\
 & &0.1&0.092&0.120 &0.103&0.969&0.936&0.110&0.117&0.814\\
 \hline
 384&512 &0.05&0.048&0.039&0.045&0.895&0.828& 0.061 & 0.065&0.699 \\
 & &0.1&0.110&0.120 &0.091&0.956&0.938&0.131&0.136&0.836\\
\hline
 384&640&0.05&0.037&0.045&0.046&0.917&0.788& 0.050 &0.062&0.702\\
 & &0.1&0.078&0.128 &0.096&0.968&0.925&0.124&0.149&0.858\\
\hline
 512&512 &0.05&0.037&0.034&0.048 &0.975&0.877& 0.044& 0.027&0.800\\
 & &0.1&0.096&0.111 &0.106&0.994&0.967&0.097&0.100&0.889\\
 \hline
 512& 640 &0.05&0.037&0.054&0.046&0.971&0.890& 0.055&0.061&0.811 \\
 & &0.1&0.094&0.137 &0.103&0.987&0.973& 0.119&0.131&0.906\\
 \hline
 640&640 &0.05&0.044&0.035&0.037 &0.993&0.959& 0.045&0.043&0.906\\
 & &0.1&0.106&0.101 &0.089&0.999&0.993&0.104&0.110&0.934\\
\hline
\end{tabular}
\caption{\textit{\label{tab1}
Rejection frequencies of the test \eqref{test} under the null hypothesis and several alternatives for $\rho=0$.}}
\end{table}

\renewcommand{\baselinestretch}{2}

\end{document}